\documentstyle[amscd,amssymb,verbatim,12pt]{amsart}
\newcommand{\lrar}[1]{\begin{picture}(50,10)(-25,-5)                          
\put(-25,0){\vector(1,0){50}}
\put(0,5){\makebox(0,0)[b]{\mbox{$#1$}}}
\end{picture}}

\newcommand{\ldar}[1]{\begin{picture}(10,50)(-5,-25)
\put(0,25){\vector(0,-1){50}}
\put(5,0){\mbox{$#1$}}
\end{picture}}

\newcommand{\ldrar}[1]{\begin{picture}(50,50)(-25,-25)
\put(-25,25){\vector(1,-1){50}}
\put(5,0){\mbox{$#1$}}
\end{picture}}

\renewcommand{\mod}{\operatorname{mod}}
\newcommand{\comod}{\operatorname{comod}}
\newcommand{\Com}{\operatorname{Com}}
\renewcommand{\Bar}{\operatorname{Bar}}
\newcommand{\tot}{\operatorname{tot}}

\newcommand{\und}{\underline}

\newcommand{\OO}{{\cal O}}
\newcommand{\Sym}{\operatorname{Sym}}

\newcommand{\DD}{{\cal D}}

\newcommand{\II}{{\cal I}}

\newcommand{\BB}{{\cal B}}

\newcommand{\mg}{{\frak m}}
\newcommand{\hra}{\hookrightarrow}
\newcommand{\lan}{\langle}
\newcommand{\ran}{\rangle}

\newcommand{\GG}{{\cal G}}
\newcommand{\CC}{{\cal C}}

\newcommand{\Spec}{\operatorname{Spec}}

\renewcommand{\P}{{\Bbb P}}

\newcommand{\Th}{\Theta}
\newcommand{\si}{\sigma}

\newcommand{\Pic}{\operatorname{Pic}}

\newcommand{\eps}{\epsilon}
\newcommand{\We}{\bigwedge}
\renewcommand{\ker}{\operatorname{ker}}

\numberwithin{equation}{subsection}

\newtheorem{thm}{Theorem}[section]
\newtheorem{prop}[thm]{Proposition}
\newtheorem{lem}[thm]{Lemma}
\newtheorem{cor}[thm]{Corollary}
\newenvironment{rem}{\vspace{3mm}\noindent
{\bf Remark.}}{\vspace{3mm}}
\newenvironment{defi}{\vspace{3mm}\noindent
{\bf Definition.}}{\vspace{3mm}}
\newenvironment{rems}{\vspace{3mm}
\noindent {\bf Remarks.}}{\vspace{3mm}}
\newenvironment{ex}{\vspace{3mm}\noindent
{\bf Example.}}{\vspace{3mm}}
\newenvironment{exs}{\vspace{3mm}\noindent
{\bf Examples.}}{\vspace{3mm}}

\newcommand{\Pf}{\noindent {\it Proof}}
\newcommand{\id}{\operatorname{id}}

\newcommand{\ov}{\overline}

\newcommand{\rk}{\operatorname{rk}}
\newcommand{\ra}{\rightarrow}
\renewcommand{\AA}{{\cal A}}
\newcommand{\FF}{{\cal F}}
\newcommand{\HH}{{\cal H}}
\newcommand{\PP}{{\cal P}}

\newcommand{\SS}{{\cal S}}
\newcommand{\LL}{{\cal L}}

\newcommand{\Hom}{\operatorname{Hom}}
\newcommand{\Ext}{\operatorname{Ext}}

\renewcommand{\a}{\alpha}
\renewcommand{\b}{\beta}
\newcommand{\om}{\omega}
\newcommand{\De}{\Delta}

\newcommand{\Z}{{\Bbb Z}}

\newcommand{\Ga}{\Gamma}

\newcommand{\wt}{\widetilde}

\newcommand{\sub}{\subset}
\newcommand{\ed}{\qed\vspace{3mm}}

\title{$A_{\infty}$-structures, Brill-Noether loci and the Fourier-Mukai
transform}
\author{A. Polishchuk}
\thanks{This research was supported in part by the NSF grant DMS-0070967}

\begin{document}
\maketitle

\bigskip

\centerline{Introduction}

\medskip

The goal of this paper is to show how the techniques of 
$A_{\infty}$-categories can be applied to the study of variation of
cohomology spaces of coherent sheaves under deformations.

Our starting point is the fact that the derived category
$D^b(X)$ of coherent sheaves on a $k$-scheme $X$,
where $k$ is a field, can be equipped with a natural structure of
$A_{\infty}$-category
(canonical up to a homotopy). The notion of $A_{\infty}$-category
generalizes the concept of $A_{\infty}$-algebra due to J.~Stasheff \cite{S}.
It was introduced by K.~Fukaya in \cite{F0} in connection with Floer
homology and then used by M.~Kontsevich in his homological formulation 
of mirror symmetry (see \cite{Kon}).
The $A_{\infty}$-structure on the derived category $D^b(X)$ can be
defined naturally using some dg-category producing $D^b(X)$
by passing to cohomology. 
This construction was first introduced by T.~V.~Kadeishvili \cite{Kad}
in the setting of $A_{\infty}$-algebras. The idea of considering
additional structures on derived categories coming from 
dg-categories goes back to \cite{BK}.

Let us assume that $X$ is projective over $k$.
The observation we make is that the $A_{\infty}$-structure 
on $D^b(X)$ can be used to describe the variation of cohomology
spaces under formal deformations of coherent sheaves on $X$.
Our main result, Theorem \ref{Fourierthm}, gives an explicit description
in terms of the $A_{\infty}$-structure
of a complex that governs such a variation over a formal neighborhood
of a given coherent sheaf in its moduli space. 
This theorem was inspired by the
study by M.~Green and R.~Lazarsfeld
of the variation of cohomology spaces under 
deformations of topologically trivial holomorphic
line bundles on a K\"ahler manifold (see \cite{GL} Thm. 3.2). 
Note that $A_{\infty}$-structures do not
appear in {\it loc.~cit.} since for topologically trivial line bundles 
on a K\"ahler manifold all higher products are homotopic to zero.

The main idea is that an $A_{\infty}$-structure on $D^b(X)$
gives rise to a canonical formal deformation of every coherent 
sheaf on $X$. For sheaves with unobstructed deformations the obtained
families are universal. In general we conjecture that they are miniversal.
A homotopy between $A_{\infty}$-structures leads to formal
changes of variables in the corresponding formal coordinate systems
on the moduli spaces. Thus, 
a choice of an $A_{\infty}$-structure can be considered as an algebraic
analogue of choosing hermitian metrics on all vector bundles,
so that the above construction is an algebraic analogue of the Kuranishi
construction (see \cite{F} for more on this analogy).
On the other hand, using the $A_{\infty}$-structure one can control
the variation of cohomology spaces in the above formal
universal families. The crucial notion that helps to organize these
deformed spaces is that of an $A_{\infty}$-functor. Namely, we show that
for every object of an $A_{\infty}$-category there is a canonical
deformation of the corresponding representable $A_{\infty}$-functor.
In the case of a coherent sheaf on $X$ this $A_{\infty}$-functor
corresponds to the canonical formal deformation mentioned above. 
Looking at the
variation of values of this $A_{\infty}$-functor on 
specific objects of $D^b(X)$
one can obtain information about the formal neighborhoods of the
loci where dimensions of cohomology jump.

As an application of our techniques we calculate
formal neighborhoods of ``sufficiently nice" points
in Brill-Noether loci parametrizing special vector bundles on curves.
Recall that the classical Brill-Noether loci for a
smooth projective curve $C$ parametrize line bundles 
of given degree on $C$ with given number
of linearly independent global sections. More precisely, for every $d\ge 0$,
$r\ge 0$, one has a subscheme $W^r_d$ in the 
Jacobian $J^d$ of degree $d$ line bundles on $C$, supported on the  
set of line bundles $L$ (of degree $d$) 
with $h^0(L)>r$ (for the precise definition see \cite{ACGH}).
Perhaps, the most important example of a Brill-Noether
locus is the theta divisor $\Th=W^0_{g-1}\sub J^{g-1}$,
where $g$ is the genus of $C$,
consisting of line bundles $L$ with $h^0(L)>0$. Riemann's theorem
asserts that the multiplicity of $\Th$ at a point $L$ is equal to
$h^0(L)$. In \cite{K} G.~Kempf generalized this theorem by
describing tangent cones to points of $W_d:=W^0_d$.
The same techniques can be used to calculate tangent cones to some
points of $W^r_d$ for $r>0$ (see \cite{ACGH}).

Similar Brill-Noether loci can be defined in the moduli
spaces of stable (or semistable) vector bundles of higher rank on $C$.
More generally, one can consider {\it twisted} Brill-Noether loci
$W^r_{n,d}(E)$ parametrizing stable vector bundles $V$ of rank $n$ and
degree $d$ such that $h^0(V\otimes E)>r$, where $E$ is a fixed vector
bundle on $C$ (see \cite{Teix}). Kempf's results admit partial
generalization to these loci (see \cite{Li} for the case $E=\OO_C$,
\cite{Teix2} for the case $n=1$, $r=0$). Using the $A_{\infty}$-techniques
we will prove the following theorem complementing these results.

\begin{thm}\label{mainthm}
Let $E$ be a vector bundle on $C$,
$V$ be a stable vector bundle on $C$ of rank $n$ and degree $d$, such that
the natural map
$$\mu_{V,E}:H^0(C,V\otimes E)\otimes H^0(C,V^{\vee}\otimes E^{\vee}\otimes\om)
\ra H^0(V\otimes V^{\vee}\otimes\om)$$
is injective. We think about $\mu_{V,E}$ as a matrix of
linear forms on $T=\Ext^1(V,V)\simeq H^0(V\otimes V^{\vee}\otimes\om)^*$.
Then the formal neighborhood of $W^r_{n,d}(E)$ 
at $V$, where $r<h:=h^0(V\otimes E)$, 
is isomorphic over $k$ to the formal neighborhood of zero in
the subscheme of $T$ defined by the 
$(h-r)\times (h-r)$ minors of $\mu_{V,E}$. 
\end{thm}

Note that the case when $n=1$, $E=\OO_C$ and $k$ is algebraically closed 
follows essentially from the definition of Brill-Noether loci (see
section \ref{BNsec}). However, already the case $n=1$, $\rk E>1$ seems
to be non-trivial.
The map $\mu_{V,E}$ is called the (generalized) {\it Gieseker-Petri map}.
Its injectivity is equivalent
to the condition that the smallest Brill-Noether locus associated with $E$
containing $V$, namely, $W^{h-1}_{n,d}(E)$, where $h=h^0(V\otimes E)$,
is smooth of expected dimension at $V$.
The above theorem describes in this
situation the formal neighborhoods of all larger Brill-Noether loci
$W^r_{n,d}(E)$, $r<h$, at $V$.
Note that the Kempf's theorem and its generalizations state that
the tangent cone to $W^r_{n,d}(\OO)$ at $V$ is isomorphic to the subscheme
of $T$ considered in Theorem \ref{mainthm} under a weaker assumption
on $V$ (see \cite{Li}).

Theorem \ref{mainthm} follows from
a stronger result, Theorem \ref{homthm}, asserting that certain
higher products associated with $V$ and $E$ are homotopic to zero.
We expect that this statement should play a role in a noncommutative
version of our results formulated in terms of canonical noncommutative
thickenings of the moduli space of vector bundles on $C$ (see \cite{Kap}).


We also apply our techniques to the study of the Fourier
transform of certain line bundles on symmetric powers $\Sym^d C$
of a curve $C$. Namely, for a line bundle $L$ on $C$ let us denote
by $L^{(d)}$ the $d$-th symmetric power of $L$ which is a line bundle
on $\Sym^d C$. Let us denote by $F_d(L)$ the derived push-forward
of $L^{(d)}$ under the natural morphism $\si^d:\Sym^d C\ra J^d$.
It is not difficult to show that if $\deg(L)\ge -1$ then 
$F_d(L)$ is actually a sheaf concentrated in degree $0$ 
(see Lemma \ref{Fdlem}(b)).  
We fix a point $p\in C$ and identify $J^d$ with $J$ by 
$L\mapsto L(-dp)$. Recall that for every abelian variety $A$
the Fourier-Mukai transform is an equivalence $\SS:D^b(A)\ra D^b(\hat{A})$,
where $\hat{A}$ is the dual abelian variety (see \cite{Mukai}).
Using the self-duality of $J$ we can consider
the Fourier-Mukai transform as an autoequivalence $\SS:D^b(J)\ra D^b(J)$.
In section \ref{symsec} we will prove the following theorem.

\begin{thm}\label{symthm} Assume that $1\le d\le g-1$.
Then one has the following isomorphisms in $D^b(J)$:
$$[-1]_J^*\SS(F_d(\OO_C((g-d)p)))\simeq F_{g-d}(\OO_C(-p))(\Th)[-d]\simeq
R\und{\Hom}(F_{g-d}(\OO_C(dp)),\OO_J),$$
where $[-1]_J:J\ra J$ is the inversion map, $\Th=W_{g-1}\sub J$ is the
theta divisor.
\end{thm}

This theorem 
provides a new collection of coherent sheaves on Jacobians
for which W.I.T. holds (see \cite{Mukai} for terminology and
for other examples).
Note that the case $d=1$ was considered in \cite{BP} in connection
with Torelli theorem.

\section{$A_{\infty}$-structures}

In this section we present some $A_{\infty}$-formalism
(for the most part, well-known).

\subsection{$A_{\infty}$-categories and functors}\label{defsec}

For more details concerning most of the following definitions the reader
can consult \cite{Keller}. 

Let $k$ be a field. All the categories (and $A_{\infty}$-categories)
considered below are going to be $k$-linear.
This means that all morphism spaces are $k$-vector spaces and all
operations are $k$-linear. By the Koszul sign rule we mean
the appearance of $(-1)^{\wt{a}\cdot\wt{b}}$ when switching
graded symbols $a$ and $b$, where we use the notation $\wt{a}=\deg(a)$.

\begin{defi} (i) An {\it $A_{\infty}$-category} 
$\CC$ consists of a class of objects
and a collection of graded morphism spaces 
$\Hom^*(O_1,O_2)=\Hom_{\CC}^*(O_1,O_2)$ for every pair
of objects $O_1$, $O_2$ equipped with the operations
$$m_n:\Hom^*(O_2,O_1)\otimes_k\Hom^*(O_3,O_2)\otimes_k\ldots
\otimes_k\Hom^*(O_{n+1},O_n)\ra\Hom^*(O_{n+1},O_1),$$
where $n=1,2,3,\ldots$, homogeneous of degree $2-n$.
These operations satisfy the following {\it $A_{\infty}$-constraint}:
$$\sum_{k+l=n+1}\sum_{j=1}^{k}
(-1)^{j+l(k-j)+\eps}
m_k(a_1,\ldots,a_{j-1},m_l(a_j,\ldots,a_{j+l-1}),a_{j+l},\ldots,a_n)=0,$$
where $n=1,2,3\ldots$, $(-1)^{\eps}$ comes from the Koszul sign rule
($m_l$ gets exchanged with $a_1,\ldots,a_{j-1}$).

\noindent (ii) An $A_{\infty}$-category $\CC$ is called {\it minimal}
if $m_1=0$.
\end{defi}

Below by a {\it non-unital} category we mean a version of the notion
of a category in which the existence of identity morphisms is
not required.

\begin{exs} 1. A (non-unital) {\it dg-category} is an 
$A_{\infty}$-category with
$m_n=0$ for $n>2$. It can be considered as a non-unital category, such that
all morphism spaces are equipped with the structure of complexes
and the composition satisfies the Leibnitz rule. On the other hand,
every minimal $A_{\infty}$-category can be considered as a non-unital
category with an additional structure given by higher products.

The most important example of a dg-category is the dg-category of complexes
$\Com(\AA)$
over some $k$-linear category $\AA$. It is defined as follows (see \cite{BK}). 
For every pair of complexes $K^{\bullet}$, $L^{\bullet}$ we set
$$\Hom^n(K^{\bullet},L^{\bullet})=
\prod_{j-i=n}\Hom_{\AA}(K^i,L^j).$$
The differential is given by $m_1(f)=d\circ f-(-1)^{\wt{f}}f\circ d$
and the composition by $m_2(f,g)=f\circ g$.

2. Let $A$ be a dg-algebra (resp. dg-coalgebra). Then we can consider
left dg-modules (resp. dg-comodules) over $A$ as objects of a dg-category.
Namely, for a pair of dg-modules (resp. dg-comodules) $M$, $M'$ we set
$\Hom^n(M,M')$ to be the space of maps $f:M\ra M'$ such that
$f(M_i)\sub M_{i+n}$ and $f$ commutes with the $A$-action (resp. coaction)
in the graded sense.
The differential $m_1$ on these spaces and the composition $m_2$
are defined by the same formulas as in the previous example.
We will denote the dg-category of dg-modules (resp. dg-comodules)
over $A$ by $A-dg-\mod$ (resp. $A-dg-\comod$). Similarly, one
can define the dg-category of right dg-modules $dg-\mod-A$
(resp. $dg-\comod-A$).
\end{exs}

\begin{defi} Let $\CC$ be an $A_{\infty}$-category. The opposite
$A_{\infty}$-category $\CC^{op}$ has the same objects as
$\CC$, the morphism spaces 
$\Hom^*_{\CC^{op}}(O_1,O_2)=\Hom^*_{\CC}(O_2,O_1)$, and the operations
$$m_n^{op}(a_1,\ldots,a_n)=(-1)^{{n+1 \choose 2}+1+\eps}m_n(a_n,\ldots,a_1),$$
where $\eps$ is determined by the Koszul sign rule:
$\eps=\sum_{i<j}\wt{a_i}\wt{a_j}$.
\end{defi}

It is easy to check that the $A_{\infty}$-constraint is indeed satisfied
for $(m_n^{op})$. The linear and the constant term in
the quadratic function defining the sign in $m_n^{op}$ are chosen in
such a way that $m_1^{op}=m_1$ and the sign $m_2^{op}$ comes only
from the Koszul sign rule.

\begin{defi} An {\it $A_{\infty}$-functor} $F:\CC\ra\CC'$ between 
$A_{\infty}$-categories associates to every object $O$ of $\CC$ an
object $F(O)$ of $\CC'$ and to every collection of objects
$O_1,\ldots,O_{n+1}$ in $\CC$, where $n\ge 1$, a $k$-linear map
$$F_n:\Hom_{\CC}^*(O_2,O_1)\otimes_k\Hom_{\CC}^*(O_3,O_2)\otimes_k\ldots
\otimes_k\Hom_{\CC}^*(O_{n+1},O_n)\ra\Hom_{\CC'}^*(F(O_{n+1}),F(O_1))$$
of degree $1-n$. These maps are compatible with the operations
in $\CC$ and $\CC'$ in the following way:
\begin{align*}
&\sum
(-1)^{\eps+d(k_{\bullet})}
m_i(F_{k_1}(a_1,\ldots,a_{k_1}),F_{k_2-k_1}(a_{k_1+1},\ldots,a_{k_2}),
\ldots,F_{n-k_{i-1}}(a_{k_{i-1}+1},\ldots,a_n))\\
&=\sum_{k+l=n+1}\sum_{j=1}^{k}(-1)^{\eps'+j-1+l(k-j)} 
F_k(a_1,\ldots,a_{j-1},m_l(a_j,\ldots,a_{j+l-1}),a_{j+l},\ldots,a_n),
\end{align*}
where in the LHS the summation is taken over all
sequences $0=k_0<k_1<k_2<\ldots<k_{i-1}<n$,
$d(k_{\bullet})=\sum_{j=1}^{i-1}(i-j)(k_j-k_{j-1}-1)$,
$\eps$ and $\eps'$ come from the Koszul sign rule.
\end{defi}

\begin{defi} Every $A_{\infty}$-functor $F:\CC\ra\CC'$ defines
the {\it opposite} $A_{\infty}$-functor $F^{op}:\CC^{op}\ra(\CC')^{op}$
between the opposite $A_{\infty}$-categories by the rule
$$F_n^{op}(a_1,\ldots,a_n)=(-1)^{{n+1 \choose 2}+1+\eps}F_n(a_n,\ldots,a_1),$$
where $\eps$ comes from the Koszul sign rule.
\end{defi}

One can define a {\it contravariant} $A_{\infty}$-functor from $\CC$ to $\CC'$
as an $A_{\infty}$-functor from $\CC$ to $(\CC')^{op}$. By the above
example this is equivalent to giving an $A_{\infty}$-functor from
$\CC^{op}$ to $\CC'$. To avoid confusion in signs we will consider only
covariant $A_{\infty}$-functors, replacing the target by the opposite 
$A_{\infty}$-category when necessary.

\begin{ex} Every object $O$ of an $A_{\infty}$-category $\CC$ defines
the {\it representable $A_{\infty}$-functor} $h_O:\CC\ra\Com(k-\mod)$,
where $\Com(k-\mod)$ is the dg-category of complexes of $k$-vector spaces.
Namely, $h_O(O')=\Hom^*(O,O')$ with the differential $m_1$ and
$$h_{O,n}(a_1,\ldots,a_n)(a)=
(-1)^{(n+1)}m_{n+1}(a_1,\ldots,a_n,a).$$
Similarly, we have the representable $A_{\infty}$-functor
$h'_A:\CC\ra\Com(k-\mod)^{op}$ defined by $h'_O(O')=(\Hom^*(O',O),m_1)$,
$$h'_{A,n}(a_1,\ldots,a_n)(a)=
(-1)^{\wt{a}(\wt{a_1}+\ldots+\wt{a_n})}
m_{n+1}(a,a_1,\ldots,a_n).$$
\end{ex}

One can define the composition of $A_{\infty}$-functors 
(see \cite{Keller} 3.4, or section \ref{barsec} below)
and the notion of a homotopy
between two $A_{\infty}$-functors with the same source and target (for
$A_{\infty}$-algebras this reduces to the notion of a
homotopy between $A_{\infty}$-morphisms defined in \cite{Keller}, 3.7; see
also \ref{barsec} below). Using this one can define the notion of
$A_{\infty}$-equivalence between $A_{\infty}$-categories.
It is known that an $A_{\infty}$-functor $F:\CC\ra\CC'$ is an 
$A_{\infty}$-equivalence
if and only if $H^*F$ is an equivalence
(see \cite{Kad2}, \cite{Pr} for the case of $A_{\infty}$-algebras).

\begin{defi} (i)
Let $\CC$ be an $A_{\infty}$-category. Then we define the 
graded non-unital category $H^*\CC$ and the non-unital category $H^0\CC$
having the same objects as $\CC$ by setting
$$\Hom_{H^*\CC}(O,O')=H^*(\Hom_{\CC}(O,O'),m_1),\
\Hom_{H^0\CC}(O,O')=H^0(\Hom_{\CC}(O,O'),m_1)$$
and by considering the composition law induced by $m_2$ for these spaces.

\noindent(ii)
The component $F_1$ of an 
$A_{\infty}$-functor $F:\CC\ra\CC'$ between $A_{\infty}$-categories
induces 
the graded non-unital functor $H^*F:H^*\CC\ra H^*\CC'$ and the non-unital
functor $H^0:H^0\CC\ra H^0\CC'$.
\end{defi}

The following theorem is essentially due to T.~V.~Kadeishvili 
(in \cite{Kad} only the case of $A_{\infty}$-algebras is considered, 
however, the generalization to $A_{\infty}$-categories is 
straightforward). It was
rediscovered several times in different contexts, see \cite{Keller}
and references therein.

\begin{thm}\label{Kadthm} 
Let $\CC$ be an $A_{\infty}$-category. 
Then there exists an extension
of the structure of graded non-unital category on $H^*\CC$ to that 
of $A_{\infty}$-category, such that $H^*\CC$ with this structure
is $A_{\infty}$-equivalent to $\CC$. More precisely, there exists
an $A_{\infty}$-functor $F:\CC\ra H^*\CC$ such that $H^*F$ is
the identity functor on $H^*\CC$. 
\end{thm}

Note that an $A_{\infty}$-structure on $H^*\CC$ constructed in the
above theorem is not unique. However, all these structures are homotopic
in the sense of the following definition.

\begin{defi} Let $\CC$ and $\CC'$ be two minimal $A_{\infty}$-categories. 
An $A_{\infty}$-functor $F:\CC\ra\CC'$ is called a
{\it homotopy} if the functor $H^*F$ is the identity. 
\end{defi}

In other words, if there is a homotopy $F:\CC\ra\CC'$ then
$\CC$ and $\CC'$ should have the same objects and the same spaces
of morphisms but possibly different sets of operations $m=(m_n)$ and 
$m'=(m'_n)$.
The fact that $F_1$ is the identity together with
the minimality assumption implies that $m_2=m'_2$, i.e. 
$\CC$ and $\CC'$ coincide as usual (non-unital) categories.
Changing the point of view, we can consider $m$ and $m'$ as two 
minimal $A_{\infty}$-structures on a non-unital category $\CC$ and 
say that $F$ is a homotopy from the $A_{\infty}$-structure $m$ to $m'$.
In fact, it is easy to see that for every $A_{\infty}$-category
$\CC$ and every collection of morphisms
$$F_n:\Hom_{\CC}^*(O_2,O_1)\otimes_k\Hom_{\CC}^*(O_3,O_2)\otimes_k\ldots
\otimes_k\Hom_{\CC}^*(O_{n+1},O_n)\ra\Hom_{\CC}^*(F(O_{n+1}),F(O_1))$$
of degree $1-n$, where $n=2,3,\ldots$, 
there exists a unique $A_{\infty}$-structure $m'=(m'_n)$ on $\CC$, such that
$F$ defines a homotopy from $m$ to $m'$ (see \cite{P-hmc}). 
We will use the notation
$m'=m+\delta(F)$ for this new $A_{\infty}$-structure.

Note that composition of two homotopies is again a homotopy.
It is easy to see that with respect
to this composition the set of all homotopies for a given non-unital
category $\CC$ forms a 
group acting on the set of all minimal $A_{\infty}$-structures on $\CC$
with $m_2$ given by the composition in $\CC$.

\subsection{$A_{\infty}$-structures on derived categories}\label{derivedsec}

Let $\AA$ be a $k$-linear abelian category with enough injective objects.
Then one can equip the derived category $D^+(\AA)$ of bounded below complexes 
with a canonical structure (up to a homotopy)
of minimal $A_{\infty}$-category such that $m_2$
is the standard composition in $D^+(\AA)$. More precisely,
first, one has to make a graded category out of $D^+(\AA)$ by taking
$$\Hom^*(K,K')=\oplus_{i\in\Z}\Hom^i_{D(\AA)}(K,K')$$
as morphism spaces, where $\Hom^i_{D(\AA)}(K,K')=\Hom_{D(\AA)}(K,K'[i])$. 
Let us denote this graded category by $D^+_{\Z}(\AA)$.
Then there is a canonical homotopy class of minimal 
$A_{\infty}$-structures
on $D^+_{\Z}(\AA)$ with $m_2$ equal to the standard composition.
It is constructed as follows. As is well-known (see \cite{GM}, III, 5.20),
the category $D^+(\AA)$ is equivalent
to the homotopy category $H^0\Com^+(\II)$
of complexes of injective objects in $\AA$ bounded below (here $\II$ denotes
the subcategory of injective objects in $\AA$).
This category is obtained by taking $H^0$ from the dg-category
$\Com^+(\II)$ of complexes.
Similarly, the $\Z$-graded category
$D^+_{\Z}(\AA)$ is equivalent to $H^*\Com^+(\II)$. Now applying 
Theorem \ref{Kadthm} we get a canonical homotopy class of 
minimal $A_{\infty}$-structures on $H^*\Com^+(\II)$.

We can restrict the above $A_{\infty}$-structure on $D^+_{\Z}(\AA)$
to the subcategory $D^b_{\Z}(\AA)$.

It is natural to ask whether the canonical homotopy class
of $A_{\infty}$-structures on $D^b_{\Z}(\AA)$ constructed above,
contains the trivial one (with $m_n=0$ for $n>2$). The
examples of computations of Massey products (see \cite{P-trisec}) show that
for derived categories of coherent sheaves on projective curves 
this is not the case---otherwise all these Massey products would vanish.
This implies that similar non-triviality holds for an arbitrary
projective variety of dimension $\ge 1$.
However, in Theorem \ref{mainthm} we assert that in some cases
the part of the $A_{\infty}$-structure on $D^b_{\Z}(C)$ (where $C$ is
a curve) responsible for the variation of cohomology near 
given stable vector bundle, is homotopically trivial in the above sense. 

\subsection{Bar construction and $A_{\infty}$-modules}\label{barsec}

Bar construction is a convenient tool to record the 
$A_{\infty}$-data. In particular, it explains the signs
arising in $A_{\infty}$-definitions and allows to define
$A_{\infty}$-morphisms and homotopies between them in a concise way.

\begin{defi} Let $A$ be an $A_{\infty}$-algebra over $k$.
Its {\it bar construction} is the space
$$\Bar(A)=T(A[1])=\oplus_{n\ge 0}(A[1])^{\otimes n}$$ 
considered as a cofree (coassociative) coalgebra (with counit)
with the coderivation $b_A:\Bar(A)\ra \Bar(A)$ of degree $1$
whose components $b_n:(A[n])^{\otimes n}\ra A[1]$, $n\ge 1$, 
are defined by
the products $m_n$ via the following commutative diagram
\begin{equation}
\begin{array}{ccc}
A^{\otimes n} &\lrar{m_n} & A\\
\ldar{s^{\otimes n}} & & \ldar{s}\\
(A[1])^{\otimes n} &\lrar{b_n} & A[1]
\end{array}
\end{equation}
where $s:A\ra A[1]$ is the canonical map of degree $-1$.
\end{defi}

Note that when the map $b_n$ applied to elements of $(A[n])^{\otimes n}$
is expressed in terms of $m_n$, some signs will arise because of the 
Koszul sign rule:
$$b_n(s(a_1)\otimes\ldots\otimes s(a_n))=
(-1)^{(n-1)\wt{a_1}+(n-2)\wt{a_2}+\ldots
+\wt{a_{n-1}}}s(m_n(a_1,\ldots,a_n)).$$
The $A_{\infty}$-constraint is equivalent to the
statement that $b_A^2=0$ (see \cite{S}), thus we can consider
$(\Bar(A),b_A)$ as a dg-coalgebra.

The importance of the bar construction ia due to the fact that
an $A_{\infty}$-morphism between $A_{\infty}$-algebras $f:A\ra A'$
is the same as a morphism of dg-coalgebras $F:\Bar(A)\ra\Bar(A')$:
one should just take the components $f_n:A^{\otimes n}\ra A'$ and make
the map $\Bar(A)\ra A'[1]$ out of them as above. This interpretation
leads to a natural definition of the composition of $A_{\infty}$-morphisms
between $A_{\infty}$-algebras.

\begin{defi} A homotopy between a pair of
$A_{\infty}$-morphisms $f,g:A\ra A'$ of $A_{\infty}$-algebras is
the morphism $H:\Bar(A)\ra\Bar(A')$ of degree $-1$ such that
$$\Delta\circ H=(F\otimes H+H\otimes G)\circ\Delta,$$
$$F-G=b\circ H+H\circ b,$$
where $F,G:\Bar(A)\ra\Bar(A')$ are morphisms of coalgebras corresponding
to $f$ and $g$.
\end{defi}

The first condition in this definition allows to recover a homotopy
$H$ from its component $\Bar(A)\ra A'$, so $H$
corresponds to a collection of maps $h_n:A^{\otimes n}\ra A'$
of degree $-n$, $n\ge 1$, satisfying some equations.
It turns out that for $A_{\infty}$-algebras over a field $k$ the 
homotopy between $A_{\infty}$-morphisms is an equivalence relation
(see \cite{Pr}). One can also consider the corresponding notion
of homotopy equivalence between $A_{\infty}$-algebras. The following
important theorem was proven by Kadeishvili (see \cite{Kad1},\cite{Kad2})
and independently by Prout\'e (see \cite{Pr}). An $A_{\infty}$-morphism 
$f=(f_n):A\ra B$ 
between $A_{\infty}$-algebras is called a {\it quasiisomorphism}
if the corresponding map $H^*f_1:H^*A\ra H^*B$ is an isomorphism.

\begin{thm} Every quasiisomorphism of $A_{\infty}$-algebras is
a homotopy equivalence.
\end{thm}

We leave to the reader to define the bar construction of an 
$A_{\infty}$-category and the notion of homotopy between 
$A_{\infty}$-functors imitating the above definitions for 
$A_{\infty}$-categories. 
The analogue of the above theorem holds also for $A_{\infty}$-categories.

\begin{defi} A (left) $A_{\infty}$-module $M$ over an $A_{\infty}$-algebra $A$
is a graded vector space over $k$ equipped with $k$-linear operations
$$m_n^M:A^{\otimes (n-1)}\otimes_k M\ra M$$
of degree $2-n$, where $n\ge 1$, satisfying the $A_{\infty}$-constraint.
Equivalently, one can say that an $A_{\infty}$-module $M$ over $A$ is
the same as an $A_{\infty}$-category $\CC_M$ 
with two objects $X, Y$ such that
$\Hom^*(Y,Y)=A$, $\Hom(X,Y)=M$, $\Hom(Y,X)=0$, $\Hom(X,X)=0$.
\end{defi}

It is easy to see that the structure of an $A_{\infty}$-module
over $A$ on a graded $k$-vector space $M$ is equivalent to the datum of
a differential $b_M$ of degree $1$ on a cofree $\Bar(A)$-comodule
$$\Bar(M):=\Bar(A)\otimes_k M[1]$$ 
such that $(\Bar(M),b_M)$ is a dg-comodule over $(\Bar(A),b_A)$, 
i.e. $b_M^2=0$ and the pair $(b_M, b_A)$
satisfies the co-Leibnitz rule. The explicit formula for
the cogenerating components of $b=b_M$ is
$$b_n(s(a_1)\otimes\ldots\otimes s(a_{n-1})\otimes s(x))=
(-1)^{(n-1)\wt{a_1}+(n-2)\wt{a_2}+\ldots+\wt{a_{n-1}}}
s(m_n(a_1,\ldots,a_{n-1},x)),
$$
where $x\in M$, $a_1,\ldots,a_{n-1}\in A$.

\begin{defi} A {\it closed morphism} of 
$A_{\infty}$-modules $f:M\ra M'$ over an
$A_{\infty}$-algebra $A$ is defined as a sequence of maps
$f_n:A^{\otimes (n-1)}\otimes_k M\ra M'$ of degree $1-n$ (where $n\ge 1$),
such that the data $(f:M\ra M',\id:A\ra A)$ defines 
an $A_{\infty}$-functor $\CC_M\ra \CC_{M'}$ identical on objects.
\end{defi}

Again we can interpret this notion in terms of the bar constructions:
an $A_{\infty}$-morphism $M\ra M'$ is the same as a closed morphism
of dg-comodules $\Bar(M)\ra\Bar(M')$ over $\Bar(A)$. Here we equip
dg-comodules over $\Bar(A)$ with the structure of a dg-category as in
section \ref{defsec}.
More generally, using the correspondence between 
$A_{\infty}$-modules
over $A$ and dg-comodules over $\Bar(A)$, we will consider 
$A_{\infty}$-modules over $A$ as objects of a dg-category denoted by
$A-\mod_{\infty}$.
By the definition, the map $M\mapsto\Bar(M)$ defines a dg-functor
$A-\mod_{\infty}\ra \Bar(A)-dg-\comod$.

Finally, let us quote the following important theorem (see \cite{Keller}, 4.2).

\begin{thm}\label{hommodthm} 
Let $f:M\ra M'$ be a closed $A_{\infty}$-morphism of
$A_{\infty}$-modules. If $f$ is a quasiisomorphism (i.e.,
induces an isomorphism $H^*(M,m_1)\ra H^*(M',m_1)$), then
$f$ is a homotopy equivalence.
\end{thm}

\section{$A_{\infty}$-structures and formal deformations}

\subsection{Completed cobar construction and
the canonical deformation of an $A_{\infty}$-module}

Let $B$ be a dg-coalgebra over $k$, $d:B\ra B$ be the corresponding 
differential of degree $1$. Then there is a natural structure of dg-algebra
on the dual graded vector space $B^*$. Namely, we define
the differential $d:B^*\ra B^*$ such that for an element $b^*\in B^*$ 
one has $(db^*)(b)=(-1)^{\wt{b}} b^*(db)$.
The multiplication on $B^*$ is given by the 
following composition
$$B^*\otimes B^*\ra (B\otimes B)^*
\stackrel{\De^*}{\ra} B^*,$$
where $\De:B\ra B\otimes B$ is the comultiplication.
Thus, for $b^*_1,b^*_2\in B^*$ one has
$$(b^*_1 b^*_2)(b)=\sum_i (-1)^{\wt{b_i}\wt{b'_i}}
b^*_1(b_i)b^*_2(b'_i),$$
where $\Delta(b)=\sum_i b_i\otimes b'_i$.
In the case $B$ has a counit $\eps:B\ra k$, the dual map $\eps^*:k\ra B^*$
will be a unit for $B^*$.

Let $A$ be an $A_{\infty}$-algebra.
Applying the above construction to the dg-coalgebra structure on
$\Bar(A)$ we obtain a dg-algebra structure on the dual space
$C(A):=\Bar(A)^*$ with the differential $c_A$ induced by $b_A$ as above.
We will call $C(A)$ the {\it completed cobar construction} of $A$.
In particular, we obtain the
associative algebra structure (with a unit) on
$$A^!=H^0(C(A),c_A)\simeq H^0(\Bar(A),b_A)^*.$$
Our notation is motivated by the non-homogeneous quadratic duality:
if $A=k\oplus A_+$ is a quadratic dg-algebra (i.e. $A$ is generated
by $A_1$ as a $k$-algebra and the defining relations are quadratic), 
then $(A_+)^!$ is
a completion of the quadratic-linear algebra dual to $A$ (see \cite{PP}). 

\begin{prop}\label{functorD}
An $A_{\infty}$-morphism $f:A\ra B$ between $A_{\infty}$-algebras
induces a homomorphism $f^!:B^!\ra A^!$ of associative algebras.
This correspondence extends to a contravariant
functor from the homotopy category
of $A_{\infty}$-algebras to the category of associative algebras.
\end{prop}

\Pf . By the definition, homotopy classes of
$A_{\infty}$-maps $f:A\ra B$ are in bijection with
homotopy classes of homomorphisms of dg-coalgebras 
$\Bar(A)\ra \Bar(B)$. It remains to use the natural functor
$C\mapsto H^0(C)^*$ from the homotopy category of coalgebras
to the category of associative algebras. 
\ed

For a pair of graded $k$-vector spaces $V$ and $W$ we set
$$\Hom_{gr}(V,W)=\oplus_{n\in\Z}\Hom(V,W)_n$$
where $\Hom(V,W)_n=\prod_{i\in\Z}\Hom(V_i,V_{i+n})$.

Let $M$ be an $A_{\infty}$-module over $A$. Then we can define a
canonical differential $c_M$ on 
$$C(M)=\Hom_{gr}(\Bar(A),M)$$ 
which makes
$(C(M),c_M)$ into a right dg-module over the dg-algebra $C(A)$
via some natural right action of $C(A)$ on $C(M)$. Let us state
this construction in slightly more general terms.

\begin{prop}\label{comodprop} 
Let $(B,d)$ be a dg-coalgebra with a counit, $M$ be a vector space.
Consider the dual dg-algebra $(B^*,d)$.
Then every dg-comodule differential
$d_M:B\otimes M\ra B\otimes M$ on a free $B$-comodule $B\otimes M$
induces naturally a dg-module differential 
$d_M^{\vee}:\Hom_{gr}(B,M)\ra\Hom_{gr}(B,M)$, where
$\Hom_{gr}(B,M)$ is considered as a right $B^*$-module via the action
$$(\phi\cdot b^*)=(\phi\otimes b^*)\circ\De,$$
where $\phi\in\Hom_{gr}(B,M)$, $b^*\in B^*$, $\phi\otimes b^*$ is
considered as map from $B\otimes B$ to $M\otimes k=M$.
This construction extends to a dg-functor from the dg-category
of dg-comodules over $B$ that are free as $B$-comodules to 
$dg-\mod-B^*$.
\end{prop}

\Pf. The differential $d_M$ is uniquely determined by its component
$\ov{d}_M:B\otimes M\ra M$. Namely, we have
$$d_M(b\otimes m)=(\id_B\otimes\ov{d}_M)(\Delta(b)\otimes m)+db\otimes m.$$
The condition $d_M^2=0$ is equivalent to $\ov{d}\circ d_M=0$, i.e.
\begin{equation}\label{beq}
\sum_i (-1)^{\wt{b_i}}\ov{d}_M(b_i\otimes\ov{d}_M(b'_i\otimes m))+
\ov{d}_M(db\otimes m),
\end{equation}
where $\De(b)=\sum_i b_i\otimes b'_i$, $m\in M$.
Now we define $d_M^{\vee}$ by the formula
$$d_M^{\vee}(\phi)(b)=\sum_i (-1)^{\wt{b_i}\wt{\phi}}
\ov{d}_M(b_i\otimes\phi(b'_i))-(-1)^{\wt{\phi}}\phi(db)$$
where $\phi\in\Hom_{gr}(B,M)$ is a homogeneous element.
It is easy to see that $c$ is a derivation of $\Hom_{gr}(B,M)$
as a right dg-module over $C(B)$. The condition
$(d_M^{\vee})^2=0$ follows easily from the equation (\ref{beq})
using coassociativity of $\Delta$. 

Every morphism of $B$-comodules $f:B\otimes M\ra B\otimes M'$
is uniquely determined by its component $\ov{f}:B\otimes M\ra M'$.
We define the corresponding morphism of right $B^*$-modules 
$f^{\vee}:\Hom_{gr}(B,M)\ra\Hom_{gr}(B,M')$ by
$$f^{\vee}(\phi)(b)=\sum_i(-1)^{\wt{\b_i}\wt{\phi}}
\ov{f}(b_i\otimes\phi(b'_i)).$$
This gives the required dg-functor.
\ed

We apply this construction to $B=\Bar(A)$ and the differential
$d_M=b_M$ on $B\otimes M$ corresponding to the $A_{\infty}$-module
structure on $M$, and call the 
obtained right dg-module $(C(M),c_M)$ over the dg-algebra $C(A)$
the {\it completed cobar construction} of $M$.

Now let us assume that an $A_{\infty}$-algebra $A$ is 
concentrated in positive degrees: $A=\oplus_{n\ge 1} A_n$.
Then $\Bar(A)$ is concentrated in nonnegative degrees, while
$C(A)$ is concentrated in nonpositive degrees.
Hence, in this case we have a surjective homomorphism of algebras
$$C(A)\ra H^0(C(A))=A^!.$$
In fact, in this case the algebra $C(A)_0$ 
$$C(A)_0=\prod_{n\ge 0}(A_1^{\otimes n})^*$$
is a completion of the tensor algebra and $A^!$ is the quotient
of $C(A)_0$ 
by the two-sided ideal generated by the image of the map
$$A_2^*\ra \prod_{n\ge 0}(A_1^{\otimes n})^*$$
with components dual to the maps 
$$(-1)^{n \choose 2}m_n:A_1^{\otimes n}\ra A_2,\ n\ge 1.$$

Let $M=\oplus_{n\in\Z}M_n$ be an $A_{\infty}$-module over $A$ 
such that all the spaces $M_n$ are finite-dimensional (in this case
we say that $M$ is {\it locally finite-dimensional}). 
We are going to define a natural $A^!$-linear differential of degree $1$ on
the free right $A^!$-module $M\otimes A^!$. 
In the case when $M$ is finite-dimensional this differential can be 
immediately obtained from the completed cobar construction of $M$.
Namely, in this case $C(M)\simeq M\otimes C(A)$, so we get a dg-module
structure on the free right $C(A)$-module $M\otimes C(A)$. 
Tensoring with $A^!=H_0(C(A))$ over $C(A)$ we obtain 
an $A^!$-linear differential on $M\otimes A^!$. In the general case 
when only graded components of $M_n$ are finite-dimensional we 
have to take the dual route.
First, we claim that the embedding $H^0\Bar(A)\sub\Bar(A)$ is
a morphism of dg-coalgebras. Let us set $B=\Bar(A)$ for brevity.
Then we have $H^0B=\ker(b:B_0\ra B_1)$, hence 
$\Delta(H^0B)$ is contained in 
$$\ker (B_0\otimes B_0\lrar{(b\otimes\id,\id\otimes b)} 
B_1\otimes B_0\oplus B_0\otimes B_1)=H^0B\otimes H^0B$$
which proves our claim. This implies that
the subspace $H^0B\otimes M\sub B\otimes M$ is preserved by
the differential $b_M$. 
Applying the construction of Proposition
\ref{comodprop} to the coalgebra $H^0B$ (with zero differential) we 
obtain the $A^!$-linear differential on 
$$\Hom_{gr}(H^0B,M)=\oplus_{n\in\Z}\Hom(H^0B,M_n)=M\otimes A^!$$
(the last equality follows from the fact that 
$M_n$ are finite-dimensional). We will denote this differential
by $c_M$ and the complex of $A^!$-modules
$(M\otimes A^!,c_M)$ by $M_{A^!}$.
Note that $A^!$ has a natural augmentation $A^!\ra k$ and
tensoring $M_{A^!}$ with $k$ over $A^!$ we obtain the
complex $(M,m_1)$. So the differential $c_M$ can be
considered as a deformation of the differential $m_1$ on $M$.

Let us write the explicit formula for the
differential $c_M$ assuming for simplicity 
that $A_1$ has countable dimension. Let $(e_1,e_2,\ldots)$ be a
basis of $A_1$, $(e_1^*,e_2^*,\ldots)$ be the dual vectors in
$A_1^*$, so that elements of $A_1^*$
are infinite series $\sum_{n=1}^{\infty} c_n e_n^*$. Then we have
\begin{equation}\label{differential}
c_M(x\otimes r)=\sum_{n\ge 0;i_1,\ldots,i_n}
(-1)^{{n+1 \choose 2}} 
m_{n+1}(e_{i_1},\ldots,e_{i_n},x)\otimes e^*_{i_1}\ldots e^*_{i_n}\cdot r 
\end{equation}
where $x\in M$, $r\in A^!$. This infinite series makes sense as an element
of $M\otimes A^!$ since $M$ is locally finite-dimensional.
Note that using this notation we can write the defining relation
in $A^!$ as follows:
\begin{equation}\label{relationR(O)}
\sum_{n\ge 1; i_1,\ldots,i_n}(-1)^{n \choose 2} 
m_n(e_{i_1},\ldots,e_{i_n})\otimes e^*_{i_1}\ldots e^*_{i_n}=0
\end{equation} 
in $A_2\otimes A^!$.

\begin{prop}\label{!functorprop}
The map $M\mapsto M_{A^!}$ extends to
a dg-functor $A-\mod^{lf}_{\infty}\ra \Com(\mod-A^!)$,
where $A-\mod^{lf}_{\infty}$ is the category of 
locally finite-dimensional $A_{\infty}$-modules.
\end{prop}

\Pf . Let $M$ and $M'$ be $A_{\infty}$-modules over $A$.
A degree $n$ morphism of $\Bar(A)$-comodules
$f:\Bar(A)\otimes M\ra \Bar(A)\otimes M'$ is determined by
its component $\ov{f}:\Bar(A)\otimes M\ra M'$. It induces a degree $n$ map
$$f^{\vee}:M\ra \Hom_{gr}(H^0\Bar(A),M')\simeq M'\otimes A^!$$
such that $f^{\vee}(m)(b)=\ov{f}(b\otimes m)$ for $m\in M$,
$b\in H^0\Bar(A)$.
We define the value of our functor on $f$ to be the
corresponding morphism of free $A^!$-modules 
$M\otimes A^!\ra M'\otimes A^!$.
It is easy to check that this is indeed a dg-functor.
\ed

\begin{defi} Let $(A,M)$ and $(A',M')$ be two pairs each consisting
of an $A_{\infty}$-algebra and an $A_{\infty}$-module over it.
We define an $A_{\infty}$-map $f:(A,M)\ra (A',M')$ as an
$A_{\infty}$-functor between the corresponding $A_{\infty}$-categories
with two objects. A homotopy between two such $A_{\infty}$-maps
is defined as a homotopy between the corresponding $A_{\infty}$-functors.
\end{defi}

Thus, we can consider the homotopy category of pairs $(A,M)$.
Let us also consider pairs of the form $(C,K)$ where $C$ is an
associative algebra, $K$ is a complex of right $C$-modules.
We define a morphism $(C,K)\ra (C',K')$ between such pairs
as a pair $(\a,\b)$, where $\a:C\ra C'$ is a
homomorphism of algebras, $\b:K'\ra K\otimes_C C'$ is a morphism
in the homotopy category of complexes of $C'$-modules. 
The Proposition \ref{functorD} can be extended to pairs $(A,M)$
as follows.

\begin{prop}\label{functorD2}
The map $(A,M)\mapsto (A^!,M_{A^!})$ extends to a contravariant functor
from the homotopy category of pairs $(A,M)$, such that $A$ 
is a positively graded $A_{\infty}$-algebra and $M$ is 
a locally finite-dimensional $A_{\infty}$-module
over it, to the category of homotopy category of
pairs consisting of an associative algebra
and a complex of right modules over it.
\end{prop}

\Pf . A homotopy class of maps $(A,M)\ra (A',M')$ 
defines a coalgebra homomorphism $H^0\Bar(A)\ra H^0\Bar(A')$ and
a homotopy class of compatible comodule morphisms 
$$H^0\Bar(A)\otimes M\ra H^0\Bar(B)\otimes M'.$$
The former map induces a homomorphism of algebras $(A')^!\ra A^!$.
We can use the component $H^0\Bar(A)\otimes M\ra M'$ of the
latter map to define a map $M\ra M'\otimes A^!$ as in the proof of
Proposition \ref{!functorprop}. The corresponding map of
free $A^!$-modules $M\otimes A^!\ra M'\otimes A^!$ is the
required chain map $M_{A^!}\ra M'_{(A')^{!}}\otimes_{(A')^!} A^!$.
\ed

\subsection{The canonical deformation of a representable $A_{\infty}$-functor}

\begin{prop} Let $\CC$ be an $A_{\infty}$-category, $O$ be an object
of $\CC$. Then $A=\Hom^*_{\CC}(O,O)$ has a natural structure of
$A_{\infty}$-algebra and the 
representable $A_{\infty}$-functor $h'_O:\CC\ra \Com(k-\mod)^{op}$
factors as the composition of an $A_{\infty}$-functor
$H'_O:\CC\ra A-\mod_{\infty}^{op}$ with the forgetting dg-functor
$A-\mod_{\infty}^{op}\ra\Com(k-\mod)^{op}$.
\end{prop}

\Pf . The structure of an $A_{\infty}$-algebra on $A$ and of an
$A_{\infty}$-module on $\Hom^*_{\CC}(X,O)$
is simply given by the operations in $\CC$. To define an $A_{\infty}$-functor
$H'_O$ extending $h'_O$ we have to define for every sequence
$(x_1:X_1\ra X_0,\ldots,x_p:X_p\ra X_{p-1})$ of morphisms in $\CC$
a morphism 
$$H'_{O,p}(x_1,\ldots,x_p):\Hom^*_{\CC}(X_0,O)\ra\Hom^*_{\CC}(X_p,O)$$
in the dg-category of $A_{\infty}$-modules over $A$.
Such a morphism corresponds to a morphism of $\Bar(A)$-comodules
$$\Bar(A)\otimes\Hom^*_{\CC}(X_0,O)[1]\ra
\Bar(A)\otimes\Hom^*_{\CC}(X_p,O)[1].$$
By the definition, this morphism is obtained by substituting $x_1,\ldots,x_p$
in the component
$$\Bar(A)\otimes\Hom^*_{\CC}(X_0,O)[1]\otimes\Hom^*_{\CC}(X_1,X_0)[1]\otimes
\ldots\otimes\Hom^*_{\CC}(X_p,X_{p-1})[1]\ra
\Bar(A)\otimes\Hom^*_{\CC}(X_p,O)[1]$$
of the differential in the bar construction of $\CC$.
Thus, the components of $H'_{O,p}(x_1,\ldots,x_p)$ have form
$$H'_{O,p}(x_1,\ldots,x_p)_n(a_1,\ldots,a_{n-1},x)=\pm
m_{n+p}(a_1,\ldots,a_{n-1},x,x_1,\ldots,x_p),$$
where $x\in\Hom^*_{\CC}(X_0,O)$, $a_1,\ldots,a_{n-1}\in A=\Hom^*_{\CC}(O,O)$.
It is not difficult to check that the axioms of an $A_{\infty}$-functor
are satisfied.
\ed

Let $O$ be an object of an $A_{\infty}$-category $\CC$. 

\begin{defi}
Let us define the associative $k$-algebra $R(O)$ (with a unit)
by setting $R(O)=(\Hom^*(O,O)_+)^!$, where
$\Hom^*(O,O)_+=\oplus_{n\ge 1}\Hom^n(O,O)$. 
\end{defi}

As we have seen above,
$R(O)$ actually depends only on $\Hom^1(O,O)$, $\Hom^2(O,O)$ and
the operations $m_n:T^n\Hom^1(O,O)\ra\Hom^2(O,O)$.

Composing the representable $A_{\infty}$-functor
$h'_O:\CC\ra Hom^*(O,O)_+-\mod_{\infty}^{op}$ with the
dg-functor $M\mapsto (M\otimes_k R(O),c_M)$ defined in Proposition
\ref{!functorprop}, we obtain the $A_{\infty}$-functor
$$F_O:\CC\ra\Com(\mod-R(O))^{op}.$$
By the definition we have
$$F_O(O')=\Hom^*(O',O)\otimes R(O)$$
with the differential given by (\ref{differential}).
The formula for the structure of an $A_{\infty}$-functor has form
$$F_{O,p}(x_1,\ldots,x_p)(x\otimes r)=
\sum_{n\ge 0;i_1,\ldots,i_n}\pm
m_{p+n+1}(e_{i_1},\ldots,e_{i_n},x,x_1,\ldots,x_p)\otimes
e^*_{i_1}\ldots e^*_{i_n}\cdot r.
$$

Note that the composition of $F_O$ with the dg-functor
$\Com(\mod-R(O))^{op}\ra\Com(\mod-k)^{op}$ given by
$M\mapsto M\otimes_{R(O)}k$ is exactly the representable
$A_{\infty}$-functor $h'_O$.
Also, the terms of all the complexes $F_O(O')$ are
free $R(O)$-modules. Thus, we can consider $F_O$ as a formal deformation
of the functor $h'_O$.

In particular, for every object $O'$ such that $\Hom^*(O',O')$ 
and $\Hom^*(O',O)$ are concentrated
in degree $0$, we get a formal deformation of the structure of right
$\Hom^0(O',O')$-module on
$\Hom^0(O',O)\otimes_k R(O)$: the deformed action of $a\in\Hom^0(O',O')$ is
given by $F_{O,1}(a)$.

Proposition \ref{functorD} easily implies
that under a homotopy of the $A_{\infty}$-structure
the algebra $R(O)$ gets replaced to an isomorphic one.
Moreover, one can check that under this isomorphism
the functor $F_O$ gets replaced by a homotopic one.

\begin{ex} Let $A$ be a complete local Noetherian commutative $k$-algebra
with the residue field $k$. Consider the derived category $D^b_{\Z}(R-\mod)$
and equip it with $A_{\infty}$-structure as in \ref{derivedsec}. 
Let $\sideset{_A}{}k$
denote $k$ considered as an $A$-module. Then there is an isomorphism
of $k$-algebras $R(\sideset{_A}{}k)\simeq A$.
\end{ex}

\subsection{Computation for derived categories}

Now we consider the derived category $\CC=D^+_{\Z}(\AA)$, where $\AA$
is an abelian category with enough injectives, so
that $\CC$ has a minimal $A_{\infty}$-structure introduced
in section \ref{derivedsec}. Let $O$ be an object of $D^+(\AA)$.
We want to compute the value of the corresponding functor
$F_O$ on an object $Q\in D^b(\AA)$ using certain adapted resolution.
Namely, let us assume that there exists a bounded above complex
$$P^{\bullet}: \ra\ldots P^{n-1}\ra P^n$$
which is quasiisomorphic to $Q$, such that for every $i\in\Z$ the space
$\Hom^*_{D(\AA)}(P^i,O)$ is concentrated in degree $0$ (for example,
below we will consider the situation where $O$ is a coherent
sheaf on a projective scheme and 
$P^i$ are sufficiently negative vector bundles). 
Then $F_O(P^i)=\Hom^0_{D(\AA)}(P^i,O)\otimes R(O)$ and we have a complex
of $R(O)$-modules
\begin{equation}\label{FAcomplex}
F_O(P^n)\ra F_O(P^{n-1})\ra\ldots
\end{equation}
with the differentials induced by the morphisms $P^i\ra P^{i+1}$.

\begin{thm}\label{quasithm} 
In the above situation the complex of $R(O)$-modules
$F_O(Q)$ is quasiisomorphic to the complex (\ref{FAcomplex}).
\end{thm}

\Pf . Let $O\ra I^{\bullet}$, $P^{\bullet}\ra J^{\bullet}$
be quasiisomorphisms with the complexes of injective objects
bounded below. We have a homotopy equivalence
of pairs
$$(\Hom_{D(\AA)}^*(O,O),\Hom_{D(\AA)}^*(Q,O))\simeq(A,\wt{M})$$
with $A:=\tot\Hom_{\AA}(I^{\bullet},I^{\bullet})$ and
$\wt{M}:=\tot\Hom_{\AA}(J^{\bullet},I^{\bullet})$, where $\tot$ denotes
the convolution of a bicomplex. 
Therefore, by Proposition \ref{functorD2} we get an isomorphism
$R(O)\simeq A^!$ and the homotopy equivalence of
complexes $F_O(Q)\simeq \wt{M}_{A^!}$. The morphism of complexes
$P^{\bullet}\ra J^{\bullet}$ induces a quasiisomorphism
of dg-modules over $A$:
$$\wt{M}\ra M:=\tot\Hom_{\AA}(P^{\bullet},I^{\bullet}).$$
It follows that $\wt{M}_{A^!}$ and $M_{A^!}$ are homotopically equivalent
as $A_{\infty}$-modules over $A$ (see Theorem \ref{hommodthm}).

Now we observe that $M_{A^!}$ is the total complex associated
with the bicomplex of $A^!$-modules
$$(\BB:=\oplus_{i,j} M^{i,j}\otimes A^!,\partial_1,\partial_2)$$
where $M^{i,j}=\Hom_{\AA}(J^{-i},I^j)$,
the differential $\partial_1$ is the $A^!$-linear map induced
by the standard map $M^{i,j}\ra M^{i+1,j}$, while the differential
$\partial_2$ is induced by the structures of dg-modules over $A$ on
the complexes $M^{i,\bullet}$. Thus, the rows of this bicomplex
are exactly the complexes $M^{i,\bullet}_{A^!}$. 
Our assumption that $\Hom^*_{D(\AA)}(P^{-i},O)=\Hom^0_{D(\AA)}(P^{-i},O)$ 
implies
that the cohomology of $M^{i,\bullet}_{A^!}$ is concentrated in degree
$0$. Therefore, the embedding of complexes of $A^!$-modules
$$(H^0(\BB,\partial_2),\partial_1)\ra \tot\BB=(\BB,\partial_1+\partial_2)$$
is a quasiisomorphism. Finally, we note that
for every $i$ there is an isomorphism
$H^0(M^{i,\bullet}_{A^!})\simeq 
\Hom_{D(\AA)}(P^{-i},O)\otimes R(O)=F(P^{-i})$ 
and the morphisms
$H^0(M^{i,\bullet}_{A^!})\ra H^0(M^{i+1,\bullet}_{A^!})$
induced by $\partial_1$ are identified with the differentials
in the complex (\ref{FAcomplex}).
\ed

\subsection{The deformation of a coherent sheaf}\label{sheafsec} 

Now let us specialize to the case when $\CC$ is the derived category
$D^b_{\Z}(X)$ of coherent sheaves on a 
a projective scheme $X$ over $k$.
Let $\FF$ be a coherent sheaf on $X$, and
let $\OO_X(1)$ be an ample line bundle on $X$.
We can consider the associative algebra $R(\FF)$ and the $A_{\infty}$-functor
$F_{\FF}:D^b_{\Z}(X)\ra\Com(\mod-R(\FF))^{op}$ constructed above.
Note that $R(\FF)$ is a quotient of the completed tensor algebra of
the space $\Ext^1(\FF,\FF)^*$. Assume that for $n>0$
the cohomology spaces $H^*(X,\OO_X(n))$ and $H^*(X,\FF(n))$
are concentrated in degree $0$. Then the formula
$$(s\otimes r)*a=F_{\FF,1}(a)(s\otimes r)=\sum_{n\ge 0;i_1,\ldots,i_n} 
(-1)^{{n\choose 2}}m_{n+2}(e_{i_1},\ldots,e_{i_n},s,a)\otimes 
e^*_{i_1}\ldots e^*_{i_n}\cdot r,$$
where $a\in H^0(X,\OO_X(m))$, $m>0$, $s\in H^0(X,\FF(m'))$, $m'>0$, 
$r\in R(\FF)$, defines the structure of a 
graded $\oplus_{n>0} H^0(X,\OO_X(n))$-module
on 
$$M_{\FF}=\oplus_{n>0} H^0(X,\FF(n))\otimes_k R(\FF).$$ 
Clearly, this structure commutes with
the $R(\FF)$-module structure on $M_{\FF}$. 

If we replace $R(\FF)$ by its abelianization $R=R(\FF)^{ab}$,
the above construction still works, so we get a structure
of a graded $\oplus_{n>0}H^0(X,\OO_X(n))\otimes R$-module on
$M=\oplus_{n>0} H^0(X,\FF(n))\otimes R.$
Then the localization $\wt{M}$ of $M$ will be a coherent
sheaf $\FF_R$ on $X\times\Spec(R)$, flat over $R$.  
Let us consider the natural homomorphism $R\ra k$ given by the
augmentation of $R$. 
The natural isomorphism
$$M\otimes_R k\simeq H^0(X,\FF(n))$$
compatible with the $\oplus_{n>0} H^0(X,\OO_X(n))$-action, induces
an isomorphism of $\FF_R|_{X\times\Spec(k)}$ with $\FF$.
Thus, the family $\FF_R$ is a deformation of $\FF$.

\begin{thm}\label{Fourierthm} 
Let $\GG\in D^b(X)$ be an object.
Then for every $n>0$ one has an isomorphism
\begin{equation}\label{mainisom}
F_{\FF}(\GG)\simeq Rp_{2*}R\und{\Hom}(p_1^*\GG,\FF_R)
\end{equation}
in the derived category of complexes of $R$-modules,
where $p_1$ and $p_2$ are the projections of the product $X\times\Spec(R)$
to its factors.
\end{thm}

\Pf . This follows easily from Theorem \ref{quasithm}. 
Indeed, let $N$ be an integer such
that $H^i(X,\FF(n))=0$ for $i>0$, $n>N$. We can choose a quasiisomorphism
$P^{\bullet}\ra\GG$, where $P^{\bullet}$ is a bounded above complex,
such that each $P^i$ is a direct sum of line bundles
$\OO_X(-n)$ with $n>N$. Then $Rp_{2*}R\und{\Hom}(p_1^*\GG,\FF_R)$
is represented by the complex of $R$-modules 
$$\ldots \ra H^0(X\times\Spec(R),(P^n)^{\vee}\otimes\FF_R)\ra
H^0(X\times\Spec(R),(P^{n-1})^{\vee}\otimes\FF_R)\ra\ldots$$ 
But this complex coincides
with the complex $\ldots\ra F_{\FF}(P^n)\ra F_{\FF}(P^{n-1})\ra\ldots$
\ed

\begin{rem} The above theorem is a generalization of
the formal analogue of Theorem 3.2 in \cite{GL}. The latter theorem
states that in the context of K\"ahler geometry the variation of cohomology
groups in a family of topologically trivial line bundles can be described
locally by a complex similar to $F_{\FF}(\OO)$ but with the differential
depending only on $m_2$. The reason for the absence of higher corrections
to this differential is that in this situation there is a natural
choice of the $A_{\infty}$-structure for which the relevant higher products
vanish (this $A_{\infty}$-structure is constructed using Dolbeault
complexes and harmonic projectors, see \cite{P-h}).
\end{rem}

\begin{cor}\label{nontrivcor}
For every $\xi\in\Ext^1(\FF,\FF)$ consider the family
$\FF_{\xi}$ over $X\times \Spec(k[\eps]/(\eps^2))$ induced by $\FF_{R}$ via
the natural homomorphism of $k$-algebras $\pi_{\xi}:R\ra k[\eps]/(\eps^2)$
defined by $\pi_{\xi}(e)=e(\xi)\cdot\eps$ for $e\in\Ext^1(\FF,\FF)^*$.
If $\xi\neq 0$ then the family $\FF_{\xi}$ is non-constant.
\end{cor}

\Pf . Let $p_1,p_2$ be the projections of the product 
$X\times \Spec(k[\eps]/(\eps^2))$ to its factors.
The above theorem implies that the object 
$$Rp_{2*}R\und{\Hom}(p_1^*\FF,\FF_{\xi})\in D^+(k[\eps]/(\eps^2)-\mod)$$ 
is represented
by the complex $\Hom^*(\FF,\FF)\otimes_k k[\eps]/(\eps^2)$
with the differential
$$d(a\otimes r)=m_2(a,\xi)\otimes r\eps,$$
where $a\in\Ext^i(\FF,\FF)$, $r\in k[\eps]/(\eps^2)$.
Thus, $d(\id_{\FF})=\xi\otimes\eps\neq 0$, hence the
$0$-th cohomology of this complex is a proper subspace of
$\Hom(\FF,\FF)\otimes_k k[\eps]/(\eps^2)$.
It follows that the dimension of $R^0p_{2*}R\und{\Hom}(p_1^*\FF,\FF_{\xi})$
over $k$ is $<2\dim_k\Hom(\FF,\FF)$. On the other hand,
for the constant family $p_1^*\FF$ we have
$$R^0p_{2*}R\und{\Hom}(p_1^*\FF,p_1^*\FF)\simeq \Hom(\FF,\FF)\otimes_k 
k[\eps]/(\eps^2)$$
which has dimension $2\dim_k\Hom(\FF,\FF)$. Hence $\FF_{\xi}$ cannot
be isomorphic to $p_1^*\FF$.
\ed

\begin{rem}
It seems plausible 
that in the above situation the family $\FF_R$ is
the miniversal formal (commutative) deformation of $\FF$. 
In the case of deformations of modules the similar statement follows
from the work of O.~A.~Laudal \cite{L}. However, it seems that
$A_{\infty}$-techniques allows to simplify calculations of {\it loc.~cit.}
We plan to return to this question and its non-commutative analogue
in a future paper. 
\end{rem}

\section{Applications}

\subsection{Brill-Noether loci}\label{BNsec}

Now let $C$ be a projective curve over a field $k$. Below we assume that
$C$ is smooth although most probably this condition can be relaxed.
Let $U(n,d)$ be the moduli space of stable bundles of rank $n$ and
degree $d$ on a curve $C$. 
Then for every vector bundle $E$ on $C$ and every $i\ge 0$
one can define a subscheme $W^r_{n,d}(E)\sub U(n,d)$ corresponding
to stable bundles $V$ such that $\dim\Hom(V\otimes F)>r$.

Consider first the case $n=1$. Then $U(1,d)=J^d$ is the Jacobian of
line bundles of degree $d$ on $C$. Let $\PP$ be the universal
family on $C\times J^d$ and let $p_1,p_2$ be the
projections of the product $C\times J^d$ to its factors. 
Then the derived push-forward
$Rp_{2*}(p_1^*E\otimes\PP)$ can be represented by a complex
$V_0\stackrel{\delta}{\ra} V_1$ of vector bundles on $J^d$. By the definition,
the ideal sheaf of the subscheme 
$W^r_{1,d}(E)\sub J^d$ is generated locally by the
$(v_0-r)\times(v_0-r)$ minors of the matrix representing $\delta$ in 
some local bases of $V_0$ and $V_1$ (here $v_0=\rk V_0$).

In the case $n>1$ the definition is similar.
The situation is complicated a little bit by the fact that 
in general there is no universal family on $C\times U(n,d)$ 
(even Zariski locally over $U(n,d)$). However, one can get around
this difficulty by working with stacks of vector bundles 
(essentially this boils down to
considering the universal family over the relevant Quot-scheme).
The reader can consult \cite{Laumon} and \cite{Mercat} for details. 

\medskip

\noindent
{\it Proof of Theorem \ref{mainthm}}.
Applying the construction of section \ref{sheafsec} to $\FF=V$
(and using the $A_{\infty}$-structure on $D^b_{\Z}(C)$), we obtain the
family
$V_R$ on $C\times \Spec(R)$, where $R=\hat{S}(\Ext^1(V,V)^*)$ is the 
completed symmetric algebra of the space $\Ext^1(V,V)^*$.
Let $\iota:\Spec(R)\ra U(n,d)$ be the corresponding morphism to the moduli
space. Then according to Corollary \ref{nontrivcor}, the tangent
map to $\iota$ at the closed point of $\Spec(R)$ is an isomorphism.
Therefore, $\iota$ induces an isomorphism of $\Spec(R)$ with the
formal neighborhood of $V$ in $U(n,d)$. 

Now applying Theorem \ref{Fourierthm} to $\FF=V$ and
$\GG=E^{\vee}$, we obtain that the object 
$Rp_{2*}(p_1^*E\otimes V_R)\in D(R-\mod)$ is represented by
the complex 
$$F_V(E^{\vee}): H^0(C,V\otimes E)\otimes_k R\stackrel{d}{\ra}
H^1(C,V\otimes E)\otimes_k R,$$
where the differential $d$ is given by (\ref{differential}).
Let us choose some bases in $H^0(C,V\otimes E)$ and
$H^1(C,V\otimes E)$ and view $d$ as an $R$-valued matrix, 
Note that the condition of injectivity of $\mu_{V,E}$ is equivalent
to surjectivity of the map
$$\Ext^1(V,V)\ra\Hom(H^0(C,V\otimes E),H^1(C,V\otimes E))$$
obtained from $\mu_{V,E}$ via Serre duality.
It follows that the leading terms of the entries of $d$ are 
linearly independent elements of $\Ext^1(V,V)^*$. Therefore,
we can choose a formal coordinate system on $U(n,d)$ at $V$,
such that the entries of $d$ will be some of the coordinate functions.
This immediately implies the result.
\ed

\begin{rems} 1. In the case when $n=1$, $E=\OO_C$ and $k$ is
algebraically closed, the assertion of
Theorem \ref{mainthm} follows easily from the fact that for every
special line bundle $L$ on $C$ there exists an effective divisor
$D$ such that the natural map $H^0(C,L)\ra H^0(C,L(D))$ is
an isomorphism and $h^1(L(D))=0$ (this trick is considered in 
details in \cite{K-ai}). Indeed, one just have to use
the resolution $\LL(D)\ra\LL(D)|_D$ for line bundles $\LL$
in a neighborhood of $L$ and apply the definition of the Brill-Noether
loci to the corresponding complex $H^0(\LL(D))\ra H^0(\LL(D)|_D)$.

\noindent
2. It seems that the condition of injectivity of $\mu_{V,E}$ in 
Theorem \ref{mainthm} can be relaxed.
For example, we checked that the conclusion
of the theorem is satisfied for double points of theta divisors in
hyperelliptic curves, even though the corresponding Petri map has
one-dimensional kernel (the details will appear elsewhere). 
\end{rems}

The moduli spaces of stable vector bundles admit canonical noncommutative
thickenings (see \cite{Kap}). Note that in section \ref{sheafsec}
we obtained naturally noncommutative deformations of coherent sheaves
and then passed to abelianization. In particular,
an $A_{\infty}$-structure gives rise to formal coordinates on the above
noncommutative thickenings of the moduli spaces of vector bundles.
We believe that there is a way to define naturally some noncommutative
thickenings of the Brill-Noether loci, so that the analogue of Theorem
\ref{mainthm} still holds for them. This should be a consequence
of the following result.

\begin{thm}\label{homthm}
Let $V$ and $E$ be vector bundles on $C$, such that
$V$ is stable and the Gieseker-Petri map
$\mu_{V,E}$ is injective. 
Then one can choose an $A_{\infty}$-structure on $D^b_{\Z}(C)$
from the canonical homotopy class in such a way that all the products
$$m_n:\Ext^1(V,V)^{\otimes (n-1)}\otimes H^0(C,V\otimes E)
\ra H^1(C,V\otimes E)$$
vanish for $n>2$.
\end{thm}

\Pf . Let us start with some $A_{\infty}$-structure on $D^b(C)$ from
the canonical homotopy class.
We want to change it to a homotopic one,
so that for the new structure the only non-zero term in
the sum defining $d$ would be the first term (involving $m_2$).
We construct the required homotopy as the infinite composition of
homotopies $(f^{(n)})$, $n=2,3,\ldots,$ where the only non-zero
component of $f^{(n)}$ is 
$$f^{(n)}_n:\Ext^1(V,V)^{\otimes n}\ra\Ext^1(V,V).$$
It is easy to see that such an infinite composition
necessarily converges. We want to choose the first map $f^{(2)}_2$ in such
a way that the following diagram would be commutative:
\begin{equation}
\begin{array}{ccc}
\Ext^1(V,V)\otimes\Ext^1(V,V) && \\
\ldar{f^{(2)}_2}& \ldrar{}& \\
\Ext^1(V,V) &\lrar{}&\Hom(H^0(V\otimes E),H^1(V\otimes E))  
\end{array}
\end{equation}
where the horizontal and the diagonal arrows are partial
dualizations of the maps given by $m_2$ and $m_3$, respectively.
In fact, by Serre duality the bottom arrow can be identified with
the dual of $\mu_{V,W}$. Hence, it is surjective, so there exists a map
$f^{(2)}_2$ making the above diagram commutative.
Let us replace the $A_{\infty}$-structure on $D^b(C)$ by 
the homotopic one: $m\mapsto m+\delta(-f^{(2)})$. 
For this new structure the map
$$m_3:\Ext^1(V,V)\otimes\Ext^1(V,V)\otimes H^0(V\otimes E)
\ra H^1(V\otimes E)$$
will be zero. Now we can choose a map $f^{(3)}_3$ which makes the following
diagram commutative: 
\begin{equation}
\begin{array}{ccc}
\Ext^1(V,V)^{\otimes 3} && \\
\ldar{f^{(3)}_3}& \ldrar{}& \\
\Ext^1(V,V) &\lrar{}&\Hom(H^0(V\otimes E),H^1(V\otimes E))
\end{array}
\end{equation}
where the horizontal arrow is the same as before, while the
diagonal arrow is given by the partial dualization of $m_4$.
Then we again replace the $A_{\infty}$-structure by the
homotopic one: $m\mapsto m+\delta(f^{(3)})$.
For this new $A_{\infty}$-structure the maps
$$m_n:\Ext^1(V,V)^{\otimes (n-1)}\otimes H^0(V\otimes E)
\ra H^1(V\otimes E)$$
will be zero for $n=3,4$. Continuing in this way we will eventually
kill all of these maps for $n\ge 3$.
\ed

\subsection{Computation of the Fourier-Mukai transform}\label{symsec}

Let $C$ be a smooth projective curve, $J^d$ the Jacobian of line bundles
of degree $d$, $\si^d:\Sym^d C\ra J^d$ the natural morphism sending
$p_1+\ldots+p_d$ to $\OO_C(p_1+\ldots+p_d)$. For a line bundle $L$
of degree $n$ we set $F_d(L)=R\si^d_*L^{(d)}$, where
$L^{(d)}$ is the $d$-th symmetric power of $L$, which is a line bundle
on $\Sym^d C$. Below we identify $J^d$ with $J$ as before using the fixed
point $p$. In particular, we consider the Brill-Noether loci $W^r_d$
as subschemes of $J$.

We need to recall some facts about the Picard group
of $\Sym^d C$. For every $d\ge 2$ there is
an exact sequence
$$0\ra\Pic(J)\ra\Pic(\Sym^d C)\stackrel{\deg}{\ra}\Z\ra 0,$$
where the embedding of $\Pic(J)$ is given by the pull-back with
respect to $\si^d$, while the map $\deg$ to $\Z$ is normalized
by the condition $\deg(L^{(d)})=\deg(L)$ for $L\in\Pic(C)$
(see \cite{Collino}, \cite{BP}). In particular, 
$\Pic^0(\Sym^d C)$ is naturally identified with $\Pic^0(J)$.
Also, for every line bundle $M$ on $\Pic(\Sym^d C)$ and
for every linear system $\P\subset\Sym^d C$ of positive
dimension, the degree of $M$ is equal to the usual degree of
$M|_{\P}$.

\begin{lem}\label{Koszullem} Let $E$ be a vector bundle on
a projective space $\P^n$ such that there is an exact sequence
$$0\ra V\otimes\OO(-1)\ra W\otimes\OO\ra E\ra 0.$$
Then $H^i(\P^n,S^j E(m))=0$ for $i>0$, $j\ge 0$, $m\ge -i$.
\end{lem}

\Pf . We have the following Koszul resolution for $S^j E$:
\begin{align*}
&0\ra\sideset{}{^j}\We V\otimes \OO(-j)\ra
\sideset{}{^{j-1}}\We V\otimes W\otimes\OO(-j+1)\ra\ldots
\ra V\otimes S^{j-1}W\otimes \OO(-1)\ra S^jW\otimes\OO\\
&\ra S^j E\ra 0.
\end{align*}
Using this resolution to compute the cohomology of $S^j E(m)$ we immediately
derive the result from the vanishing of 
$H^i(\P^n,\OO(m))$ for $i>0$, $m\ge -i$.
\ed

\begin{lem}\label{Fdlem} (a) For every line bundle $\xi$ of degree $0$
on $C$ one has 
$$F_d(L\otimes\xi)\simeq F_d(L)\otimes \PP^{-1}_{\xi}$$
where $\PP_{\xi}$ is the line bundle on $J$ corresponding to $\xi\in J$
via the self-duality of $J$.

\noindent
(b) If $\deg(L)\ge -1$ then $F_d(L)$ is a sheaf on $J^d$,
i.e., $R^i\si^d_*L^{(d)}=0$ for $i>0$.

\noindent
(c) Assume that $d\le g+1$.
Let $j:J\setminus W^1_d\hra J$ be the open embedding. If
$-1\le \deg(L)\le g-d$ then $F_d(L)\simeq j_*j^*F_d(L)$. If 
$C$ is not hyperelliptic then this is true for every $L$ such that
$\deg(L)\ge -1$.
\end{lem}

\Pf . (a) It suffices to prove that 
$$(L\otimes\xi)^{(d)}\simeq L^{(d)}\otimes(\si^d)^*\PP^{-1}_{\xi}.$$
Since the line bundle $(L\otimes\xi)^{(d)}\otimes (L^{(d)})^{-1}$ 
on $\Sym^d C$ is algebraically equivalent to the trivial bundle, it
has form $(\si^d)^*\LL$ for some $\LL\in\Pic^0(J)$. Let
us embed $C$ into $\Sym^d C$ by the map $x\mapsto x+(d-1)p$. Restricting
our line bundles to $C$ we obtain an isomorphism 
$\xi\simeq(\si^d)^*\LL$, which precisely means that $\LL\simeq\PP^{-1}_{\xi}$.
Note that the sign appears here because the pull-back map
$\Pic^0(J)\ra\Pic^0(C)\simeq J(k)$ differs from the standard isomorphism
of $\hat{J}$ with $J$ by $[-1]_J$.

\noindent
(b) Recall that the fibers of $\si^d$ are projective
spaces corresponding to complete linear systems of degree $d$.
Note that the restriction of $L^{(d)}$ to every
fiber $\P=(\si^d)^{-1}(\xi)$, where $\xi\in J$, has degree $\deg(L)$. 
Also, it is well-known that the normal bundle $N$ to
the embedding $\P=(\si^d)^{-1}(\xi)\sub\Sym^d C$ fits into the exact sequence
$$0\ra N\ra H^1(C,\OO_C)\otimes\OO_{\P}\ra 
H^1(C,\xi(dp))\otimes\OO_{\P}(1)\ra 0.
$$
Using Lemma \ref{Koszullem} we deduce that higher cohomology of the bundle
$L^{(d)}|_{\P}\otimes S^j(N^{\vee})$ on the projective space $\P$ vanish.
By the formal functions theorem this implies that
$R^i\si^d_*L^{(d)}|_{\xi}=0$ for $i>0$. 

\noindent
(c) Let $Q_d\subset \Sym^d$ be the preimage of $W^1_d$ under the morphism
$\si^d$. Let us denote by $j'$ the open embedding 
$\Sym^d C\setminus Q_d\hra\Sym^d C$. Then it suffices to prove that
the natural map
$$\si^d_*L^{(d)}\ra \si^d_*j'_*j^{\prime *}L^{(d)}$$
is an isomorphism.
If $C$ is non-hyperelliptic then $Q_d$ has codimension $\ge 2$ in
$\Sym^d C$ (this follows from Martens theorem, see \cite{ACGH}),
so in this case $L^{(d)}\simeq j'_*j^{\prime *}L^{(d)}$.
Now assume that $C$ is hyperelliptic. Then $Q_d$ is an irreducible divisor
in $\Sym^d C$. It is easy to check that $\deg(Q_d)=d-g-1$ (see the proof
of Lemma 2.5 in \cite{BP}). Therefore, for every $n>0$ we have
$\deg(L^{(d)}(nQ_d))<0$. 
Let us consider the open subset $U=\Sym^d C\setminus Q'_d\sub \Sym^d C$, 
where $Q'_d=(\si^d)^{-1}(W^2_d)$. Since
the morphism $\pi:Q_d\cap U\ra W^1_d\setminus W^2_d$ induced
by $\si^d$ is flat, the base change theorem implies that 
$$\pi_*(L^{(d)}(nQ_d)|_{Q_d\cap U})=0$$
for every $n>0$. Hence, the natural map
\begin{equation}\label{sidU}
\si^d_{U*}(L^{(d)}|_U)\ra\si^d_{U*}(L^{(d)}(*Q_d)|_U),
\end{equation}
is an isomorphism,
where $\si^d_U:U\ra W_d\setminus W^2_d$ is the map induced by $\si^d$,
$L^{(d)}(*Q_d)=\cup_n L^{(d)}(nQ_d)$. Let $j_U:U\hra\Sym^d C$ and
$j'':\Sym^d C\setminus Q_d\ra U$ be the open embeddings, so that
$j'=j_U\circ j''$. Since the codimension of $Q'_d$ in $\Sym^d C$ is
$\ge 2$, we have $L^{(d)}\simeq j_{U*}(L^{(d)}|_U)$ 
(resp. $L^{(d)}(*Q_d)\simeq j_{U*}(L^{(d)}(*Q_d)|_U)$).
Therefore, applying the push-forward with respect to the embedding
$W_d\setminus W^2_d\hra W_d$ to the isomorphism (\ref{sidU}) we obtain the 
isomorphism
$$\si^d_*L^{(d)}\wt{\ra}\si^d_*L^{(d)}(*Q_d).$$
\ed

Note that the morphism $\si^d$ induces an isomorphism
$$H^1(J,\OO_J)\simeq H^1(\Sym^d C,\OO_{\Sym^d C}).$$
Hence, both spaces are naturally isomorphic to $H^1(C,\OO_C)$.
This isomorphism is used in the following lemma.

\begin{prop}\label{cohlem} 
Let $L$ be a line bundle on $C$. For every
$d\ge 1$, $i\ge 0$, there is a canonical isomorphism
$$H^i(\Sym^d C,L^{(d)})\simeq TS^{d-i} H^0(C,L)\otimes\sideset{}{^i}\We 
H^1(C,L),$$
where for every vector space $V$ we denote by $TS^n V\subset V^{\otimes n}$
the space of symmetric $n$-tensors. 
Under these isomorphisms the cup-product maps
$$H^1(\Sym^d C,\OO_{\Sym^d C})\otimes 
H^i(\Sym^d C,L^{(d)})\ra H^{i+1}(\Sym^d C,L^{(d)})$$ 
get identified with the natural maps
$$H^1(\OO_C)\otimes 
TS^{d-i} H^0(L)\otimes\sideset{}{^i}\We H^1(L)\ra
TS^{d-i-1} H^0(L)\otimes\sideset{}{^{i+1}}\We H^1(L)$$
induced by the cup-product map $H^1(\OO_C)\otimes H^0(L)\ra H^1(L)$.  
\end{prop}

\Pf . The first assertion is a consequence of the
{\it symmetric K\"unneth isomorphism} constructed by
Deligne in \cite{SGA4}, XVII.(5.5.17.2), (5.5.32.1). This is
a canonical isomorphism of graded vector spaces 
\begin{equation}\label{arrow}
TS^d R\Ga(C,L)\ra R\Ga(\Sym^d C, L^{(d)})
\end{equation}
where in the LHS we take the $d$-th symmetric power of the graded
vector space $R\Ga(C,L)=H^0(C,L)\oplus H^1(C,L)[-1]$, so that
$$TS^d R\Ga(C,L)\simeq\oplus_{i} TS^{d-i} H^0(C,L)\otimes\sideset{}{^i}\We 
H^1(C,L).$$
The compatibility with cup-products can be easily checked by considering
pull-backs to the $d$-th cartesian power of $C$.
\ed

\begin{cor}\label{cohcor} Assume that $\deg(L)\ge g-d$ and $h^0(L)=0$.
Then $H^i(\Sym^d C,L^{(d)})=0$ for all $i$.
\end{cor}

Let $\Th=W_{g-1}\sub J$ be the theta divisor.
The following lemma is probably well-known, however, we could
not find the reference in the literature (the case $d=1$ 
is the classical Riemann's theorem
on intersection of $C$ with $\Th$).

\begin{lem}\label{canbunlem} For every $d\ge 1$ one has an isomorphism
$$\om_{\Sym^dC}\simeq(\si^d)^*(\OO_J(\Th))((g-d-1)R_p^d),$$
where $R_p^d\subset\Sym^d C$ is the image of the natural embedding
$\Sym^{d-1} C\ra\Sym^d C:D\mapsto p+D$.
\end{lem}

\Pf . Assume first that $d>2g-2$. Then $\si^d$ identifies $\Sym^dC$
with the projective bundle $\P(E_d)$ of lines in a vector bundle
$E_d$ over $J$ in such a way that the line bundle $\OO_{\Sym^d C}(R_p^d)$
corresponds to $\OO(1)$. Furthermore, it is known
that $\det E_d\simeq\OO_J(-\Th)$ (see \cite{Mat}). This easily implies
the result in this case. To deduce it for all $d\ge 1$ one can use
the descending induction together with the fact that the normal bundle
to the embedding $\Sym^dC\ra\Sym^{d+1}C:D\mapsto p+D$ is isomorphic
to $\OO_C(R_p^d)$.
\ed 

\begin{lem}\label{locallem} 
Let $A$ be a local regular ring,
$P_0\ra P_1\ra\ldots P_n$ be a complex of free
$A$-modules of finite rank. Assume that all modules
$H^i(P_{\bullet})$ have support of codimension $\ge n$.
Then $H^i(P_{\bullet})=0$ for $i<n$.
\end{lem}

\Pf . Let $g$ be the dimension of $A$.
First of all, note that the cases $n=1$ and $n>g$ are trivial.
On the other hand,
the case $n=g$ is equivalent to the lemma in \cite{Mumford}, III.13.
To deduce the general case we will use the induction in $g$.
We can assume that $n<g$. 
Let us choose an element $x$ in the maximal ideal $\mg\subset A$,
such that $x$ does not belong to $\mg^2$ and to all associate primes
of height $n$ of the modules $(H^i(P_{\bullet}))$ (such an element
exists since $n<g$). Then
the ring $\ov{A}=A/xA$ is regular of dimension $g-1$ and we claim that the
assumptions of lemma are satisfied for the complex 
$\ov{P}_{\bullet}=P_{\bullet}/xP_{\bullet}$ of free $\ov{A}$-modules.
Indeed, from the long exact sequence
$$\ldots\ra 
H^i(P_{\bullet})\stackrel{x}{\ra} 
H^i(P_{\bullet})\ra H^i(\ov{P}_{\bullet})\ra
H^{i+1}(P_{\bullet})\stackrel{x}{\ra} 
H^{i+1}(P_{\bullet})\ra\ldots$$
we immediately obtain that the modules $H^i(\ov{P}_{\bullet})$ have
codimension $\ge n+1$ in $\Spec(A)$, or equivalently,
codimension $\ge n$ in $\Spec(\ov{A})$.
By induction assumption $H^i(\ov{P}_{\bullet})=0$ for $i<n$. Now
the above exact sequence implies that the endomorphism
of multiplication by $x$ on $H^i(P_{\bullet})$ is surjective for $i<n$.
By Nakayama lemma this implies that $H^i(P_{\bullet})$ vanishes
for $i<n$.  
\ed 

\noindent
{\it Proof of Theorem \ref{symthm}}.
The fact that $F_d(\OO_C((g-d)p))$ is concentrated in degree $0$
was already proven in lemma \ref{Fdlem}(b). Let us denote for brevity
the functor $[-1]^*\circ\SS$ by $\SS^-$.
Note that the transform $\SS^-(F_d(\OO_C((g-d)p)))$ is isomorphic
to 
$$Rp_2^*(p_1^*(\OO_C((g-d)p))^{(d)}\otimes (\sigma^d\times\id)^*\PP^{-1}),$$
where $p_1$ and $p_2$ are projections of the product
$\Sym^d C\times J$ to its factors, $\PP$ is the Poincar\'e line bundle
on $J\times J$. It follows that $\SS^-(F_d(\OO_C((g-d)p)))$ 
can be represented locally on $J$ by a complex of vector bundles
\begin{equation}\label{complexV}
V_0\ra V_1\ra\ldots\ra V_d.
\end{equation}
Now Lemma \ref{Fdlem}(a) and Corollary \ref{cohcor} imply that the 
cohomology sheaves of $\SS^-(F_d(\OO_C((g-d)p)))$ are supported on 
$W_{g-d}=\si^{g-d}(\Sym^d C)\sub J$.
Applying Lemma \ref{locallem} we derive that the cohomology of
the complex $V_{\bullet}$ is concentrated in degree $d$. Thus, 
$$S_d:=\SS^-(F_d(\OO_C((g-d)p)))[d]$$ 
is actually a coherent sheaf (placed in degree $0$).

We claim that the sheaf
$S_d$ is obtained as the (non-derived) push-forward of its
restriction to the open subset $J\setminus W^1_{g-d}$ 
(since $W^1_1=\emptyset$, we can assume below that $g-d\ge 2$).
Indeed, since locally we have a resolution of $S_d$ by
vector bundles (\ref{complexV}) of length $d$, it follows that
$\HH^i_Y(S_d)=0$ for every closed subset $Y\subset J$ of codimension
$>i+d$ (where $\HH^i_Y$ denotes local cohomology with support at $Y$).
Since $\dim W^1_{g-d}\le g-d-2$, we obtain that
$\HH^i_{W^1_{g-d}}(S_d)=0$ for $i=0,1$. This implies our claim that
\begin{equation}\label{Sdext}
S_d\simeq j_* j^*S_d, 
\end{equation}
where $j:J\setminus W^1_{g-d}\ra J$
is the natural embedding.

On the other hand, we can use 
Theorem \ref{Fourierthm} and Lemma \ref{cohlem}
to study the sheaf $S_d$ near a point $L(-(g-d)p)$ of 
$W_{g-d}\setminus W^1_{g-d}\sub J$, where $\deg(L)=g-d$ 
(so that $h^0(L)=1$ and $h^1(L)=d$).
More precisely, choosing a non-zero section $s$ of $L$,
we obtain that in a formal neighborhood of this point
$S_d$ is isomorphic to the $d$-th cohomology of the complex
$$R\ra H^1(L)\otimes R\ra\ldots
\ra \sideset{}{^{d-1}}\We H^1(L)\otimes R
\ra \sideset{}{^d}\We H^1(L)\otimes R,$$
where $R=\hat{S}(H^1(\OO)^*)$.
The differential is the sum of the usual Koszul differential associated with
the surjection $H^1(\OO_C)\ra H^1(L)$ (given by the cup-product
with $s$) and some terms that vanish modulo $\mg^2$, where $\mg\sub R$
is a maximal ideal. Therefore, its $d$-th cohomology is isomorphic
to $R/(f_1,\ldots,f_d)$ with $f_i\in\mg$, where $(f_1,\ldots,f_d)\mod\mg^2$
are components of the above map $H^1(\OO)\otimes H^1(L)$ with respect
to some basis of $H^1(L)$. It follows that the ideal $(f_1,\ldots,f_d)$
is simple and has height $d$. Since $S_d$ is supported on $W_{g-d}$
which is smooth at $L(-(g-d)p)$ and also has codimension $d$, the ideal
$(f_1,\ldots,f_d)\sub R$ coincides with the ideal defining $W_{g-d}$ near
$L(-(g-d)p)$. Therefore, $j^*S_d$ is actually a line bundle on 
$W_{g-d}\setminus W^1_{g-d}$.

Next, we want to study the derived pull-back $L(\si^{g-d})^*S_d$.
First, we need to calculate the pull-back of the Poincar\'e line
bundle $\PP$ on $J\times J$ under the morphism
$$\si^d\times\si^{g-d}:\Sym^d C\times\Sym^{g-d} C\ra J\times J.$$
For this purpose it is convenient to use the Deligne symbol
of a pair of line bundles on a relative curve (see \cite{Deligne}).
Namely, it is well-known that 
$$\PP^{-1}\simeq \lan p_{13}^*\PP_C, p_{23}^*\PP_C\ran$$
where $p_{ij}$ are projections from the product $C\times J\times J$,
$\PP_C$ is the Poincar\'e line bundle on $C\times J$ (which we
always take to be normalized at $p$). Therefore,
\begin{align*}
&(\si^d\times\si^{g-d})^*\PP^{-1}\simeq\\
&\lan \OO(\DD_{12}-d[p\times \Sym^d C\times\Sym^{g-d}C]),
\OO(\DD_{13}-(g-d)[p\times \Sym^d C\times\Sym^{g-d}C])\ran,
\end{align*}
where we consider $C\times \Sym^d C\times\Sym^{g-d}C$ as
a relative curve over $\Sym^d C\times\Sym^{g-d}C$,
$$\DD_{12}=\{(x,D,D')\in C\times\Sym^d C\times\Sym^{g-d}C: x\in D\},$$ 
$$\DD_{13}=\{(x,D,D')\in C\times\Sym^d C\times\Sym^{g-d}C: x\in D'\}.$$
Note that the intersection $\DD_{12}\cap\DD_{13}$ is irreducible
and projects birationally to the divisor $\DD\sub\Sym^d C\times\Sym^{g-d} C$
supported on the set of $(D,D')$ such that $D\cap D'\neq\emptyset$.
It follows that
$$\lan \OO(\DD_{12}),\OO(\DD_{13})\ran\simeq
\OO_{\Sym^d C\times\Sym^{g-d} C}(\DD).$$
Similarly, we derive that
$$\lan \OO(\DD_{12}),\OO(p\times\Sym^d C\times\Sym^{g-d} C)\ran
\simeq\OO_{\Sym^d C\times\Sym^{g-d} C}(R^d_p\times\Sym^{g-d}C),$$
$$\lan \OO(\DD_{13}),\OO(p\times\Sym^d C\times\Sym^{g-d} C)\ran
\simeq\OO_{\Sym^d C\times\Sym^{g-d} C}(\Sym^d C\times R^{g-d}_p),$$
where we use the notation of Lemma \ref{canbunlem}.
Thus, we obtain
$$(\si^d\times\si^{g-d})^*\PP^{-1}\simeq
\OO_{\Sym^d C\times\Sym^{g-d} C}(\DD-(g-d)[R_p^d\times\Sym^{g-d} C]-
d[\Sym^d C\times R_p^{g-d}]).
$$
Note that $\OO_{\Sym^i C}(R_p^i)$ is isomorphic
to the $i$-th symmetric power of $\OO_C(p)$.
Therefore, we derive that
\begin{align*}
&L(\si^{g-d})^*\SS^-(F_d(\OO_C((g-d)p)))\simeq
Rp_{2*}((\sigma^d\times\si^{g-d})^*\PP^{-1}((g-d)[R_p^d\times\Sym^{g-d} C]))
\simeq\\
&Rp_{2*}(\OO(\DD))(-dR_p^{g-d}),
\end{align*}
where $p_1, p_2$ are the projections of the product
$\Sym^d C\times\Sym^{g-d} C$ onto its factors.
It is easy to see that for every $D\in\Sym^{g-d}C$, the restriction of
the divisor $\DD$ to $\Sym^dC\times\{D\}$ is the divisor 
$R_D^d\sub\Sym^d C$ such that $\OO_{\Sym^dC}(R_D^p)\simeq (\OO_C(D))^{(d)}$.
It follows that over the complement to $Q^{g-d}\subset\Sym^{g-d}(C)$
the natural map $\OO_{\Sym^{g-d}C}\ra Rp_{2*}(\OO(\DD))$
is an isomorphism (by the base change theorem).
Therefore, we obtain
$$Li^*j^*S_d[-d]\simeq 
j^*\si^{g-d}_*(\OO_{\Sym^{g-d}C}(-dR_p^{g-d})),$$
where 
$i:W_{g-d}^{ns}=W_{g-d}\setminus W^1_{g-d}\ra J\setminus W^1_{g-d}$ 
is the closed embedding,
$j:J\setminus W^1_{g-d}\ra J$ is the open embedding.
On the other hand, applying the duality theory we get
$$Li^*j^*S_d[-d]\simeq j^*S_d\otimes\om^{-1}_{W_{g-d}^{ns}}.$$
Hence, the above isomorphism can be rewritten as
$$j^*S_d\simeq j^*\si^{g-d}_*(\om_{\Sym^{g-d}C}(-dR_p^{g-d})).$$
Now Lemma \ref{Fdlem}(c) together with (\ref{Sdext}) imply that
$$S_d\simeq\si^{g-d}_*(\om_{\Sym^{g-d}C})(-dR_p^{g-d}).$$
Using Lemma \ref{canbunlem} we obtain that
$$S_d\simeq\si^{g-d}_*(\OO_{\Sym^{g-d}C})(-\Th)$$
which proves the first isomorphism of the theorem.
Finally, the duality for the morphism $\si^{g-d}$ implies
that 
$$R\si^{g-d}_*(\om_{\Sym^{g-d}C}(-dR_p^{g-d}))[-d]\simeq
R\und{\Hom}(R\si^{g-d}_*\OO_{\Sym^{g-d}C}(dR_p^{g-d}),\OO_J)$$
which gives the second isomorphism of the theorem.
\ed


\end{document}